\documentclass[12pt,a4paper]{article}

\usepackage[english]{babel}
\usepackage{amssymb}
\usepackage[latin1]{inputenc}
\usepackage{graphicx}
\usepackage{amsthm}
\usepackage{amsmath}
\usepackage{color}
\usepackage{stmaryrd}
\usepackage{enumerate}

\date{}

\begin{document}

\title{On convergence criteria for incompressible Navier-Stokes equations with Navier boundary conditions and physical slip rates}

\maketitle

\begin{center}

Yasunori Maekawa\\
Department of Mathematics, Graduate School of Science, Kyoto University\\
maekawa@math.kyoto-u.ac.jp

\vspace{0.3cm}

Matthew Paddick\\
Sorbonne Universit\'es, UPMC Univ Paris 06, UMR 7598, Laboratoire Jacques-Louis Lions\\
paddick@ljll.math.upmc.fr

\end{center}

\renewcommand{\labelitemi}{$\bullet$}
\newtheorem{theo}{Theorem}
\newtheorem{lemma}[theo]{Lemma}
\newtheorem{propo}[theo]{Proposition}
\newtheorem*{defi}{Definition}
\newtheorem*{nota}{Notations}
\newtheorem{cor}{Corollary}

\newcommand{\derp}[1]{\partial_{#1}}
\renewcommand{\div}{\mathrm{div}~}
\newcommand{\eps}{\varepsilon}
\newcommand{\norme}[2]{\left\| #2 \right\| _{#1}}
\newcommand{\rot}{\mathrm{rot}~}

\maketitle

\begin{abstract}
We prove some criteria for the convergence of weak solutions of the 2D incompressible Navier-Stokes equations with Navier slip boundary conditions to a strong solution of incompressible Euler.
 The slip rate depends on a power of the Reynolds number, and it is increasingly apparent that the power 1 may be critical for $L^2$ convergence, as hinted at in \cite{Pm11}.
\end{abstract}

\section{The inviscid limit problem with Navier-slip boundary conditions}

In this brief note, we shed some light on how some well-known criteria for $L^2$ convergence in the inviscid limit for incompressible fluids work when the boundary condition is changed.
 We consider the two-dimensional Navier-Stokes equation on the half-plane $\Omega = \{(x,y)\in\mathbb{R}^2~|~ y>0\}$,
 \begin{equation} \left\{ \begin{array}{rcl} \derp{t}u^\eps + u^\eps\cdot\nabla u^\eps - \eps \Delta u^\eps + \nabla p^\eps & = & 0 \\ \div u^\eps & = & 0 \\ u^\eps|_{t=0} & = & u^\eps_0, \end{array} \right. \label{NS} \end{equation}
and study the inviscid limit problem. This involves taking $\eps\rightarrow 0$, and the question of whether the solutions of (\ref{NS}) converge towards a solution of the formal limit, the Euler equation,
 \begin{equation} \left\{ \begin{array}{rcl} \derp{t}v + v\cdot\nabla v + \nabla q & = & 0 \\ \div v & = & 0 \\ v|_{t=0} & = & v_0, \end{array} \right. \label{IE} \end{equation}
 in presence of a boundary is one of the most challenging in fluid dynamics. This is because the boundary conditions required for (\ref{IE}) are different to those for (\ref{NS}).
 In the inviscid model, there only remains the non-penetration condition
 \begin{equation} v\cdot n|_{y=0} = v_2|_{y=0} = 0, \label{EBC} \end{equation}
 hence inviscid fluids are allowed to slip freely along the boundary, while viscous fluids adhere to it when the most commonly used boundary condition, homogeneous Dirichlet,
 \begin{equation} u^\eps|_{y=0} = 0, \label{DBC} \end{equation}
 is used. As $\eps$ goes to zero, solutions of the Navier-Stokes equation are expected to satisfy the following ansatz,
 $$ u^\eps(t,x,y) = v(t,x,y) + V^\eps\left(t,x,\frac{y}{\sqrt\eps}\right), $$
 where $V^\eps$ is a boundary layer, such that $V^\eps(t,x,0) = -v(t,x,0)$.

However, the validity of such an expansion is hard to prove, and, in some cases, such as when $v$ is a linearly unstable 1D shear flow, it is wrong in the Sobolev space $H^1$, as shown by E. Grenier \cite{Ge}.
 General validity results require considerable regularity on the data. M. Sammartino and R. Caflisch proved the stability of Prandtl boundary layers in the analytic case \cite{SC2}, and the first author \cite{Mae} proved in the case when the initial Euler vorticity is located away from the boundary.
Recently, this has been extended to Gevrey framework by the first author in collaboration with D. G\'erard-Varet and N. Masmoudi \cite{GVMM}.
Precisely, in \cite{GVMM} a Gevrey stability of shear boundary layer  is proved  when the shear boundary layer profile satisfies some monotonicity and concavity conditions. One of the main objectives there is the system  
\begin{equation}\label{eq.gvmm} \left\{ \begin{array}{rcl} \derp{t}v^\eps  - \eps \Delta v^\eps + V^\eps \partial_x v^\eps + v^\eps_2 \partial_y V^\eps {\bf e}_1+ \nabla p^\eps & = & -v^\eps\cdot\nabla v^\eps\,,  \\ \div v^\eps & = & 0\,, \\ v^\eps|_{y=0} =0\,, \qquad v^\eps|_{t=0} & = & v^\eps_0\,. \end{array} \right.  \end{equation}
Here $V^\eps (y) = U^E (y) -U^E(0) + U (\frac{y}{\sqrt{\eps}})$, and $(U^E,0)$ describes the outer shear flow and $U$ is a given boundary layer profile of shear type.
In \cite{GVMM} the data is assumed to be periodic in $x$, and the following Gevrey class is introduced:
\begin{align}
\begin{split}
X_{\gamma,K} & \, = \, \{ f\in L^2_\sigma (\mathbb{T} \times \mathbb{R}_+)~|~\\
&  \| f\|_{X_{\gamma,K}} = \sup_{n\in \mathbb{Z}} (1+|n|)^{10} e^{K|n|^\gamma} \| \hat{f} (n,\cdot)\|_{L^2_y (\mathbb{R}_+)}<\infty \}\,.
\end{split}
\end{align}
Here $K>0$, $\gamma\geq 0$, and $\hat{f}(n,y)$ is the $n$th Fourier mode of $f (\cdot,y)$.  
The key concavity condition on $U$ and the regularity conditions on $U^E$ and $U$ are stated as follows:

\vspace{0.3cm}

{\bf (A1)} $U^E, U\in BC^2(\mathbb{R}_+)$, and $\displaystyle \sum_{k=0,1,2} \sup_{Y\geq 0} (1+Y^k) | \partial_Y^k U (Y)| <\infty$. 

\vspace{0.1cm}

{\bf (A2)} $\partial_Y U>0$ for $Y\geq 0$, $U(0)=0$, and $\displaystyle \lim_{Y\rightarrow \infty} U(Y) = U^E(0)$.

\vspace{0.1cm}

{\bf (A3)} There is $M>0$ such that \, $-M\partial_Y^2 U\geq (\partial_Y U)^2$ for $Y\geq 0$.

\begin{theo}[\cite{GVMM}]\label{thm.gvmm} Assume that {\bf (A1)}-{\bf (A3)} hold. Let $K>0$, $\gamma\in [\frac23,1]$. Then there exist $C, T, K', N>0$ such that for all small $\eps$ and $v_0^\eps\in X_{\gamma,K}$ with $\| v_0^\eps\|_{X_{\gamma,K}}\leq \eps^N$, the system \eqref{eq.gvmm} admits a unique solution $v^\eps\in C([0,T]; L^2_\sigma (\mathbb{T}\times \mathbb{R}_+))$ satisfying the estimate
\begin{align}
\sup_{0\leq t\leq T} \bigg ( \| v^\eps (t) \|_{X_{\gamma,K'}} + (\eps t)^\frac14 \| v^\eps (t) \|_{L^\infty} + (\eps t)^\frac12 \| \nabla v^\eps (t) \|_{L^2} \bigg ) \leq C \| v_0^\eps \|_{X_{\gamma,K}}\,.
\end{align}
\end{theo}
In Theorem \ref{thm.gvmm} the condition $\gamma\geq \frac23$ is optimal at least in the linear level, due to the Tollmien-Schlichting instability; see Grenier, Guo, and Nguyen \cite{GGN2014}. More general results, including the case when $U^E$ and $U$ depend also on the time variable, can be obtained; see \cite{GVMM} for details. 
\vspace{12pt}

The situation remains delicate when the Dirichlet boundary condition (\ref{DBC}) is replaced by (\ref{EBC}) plus a mixed boundary condition such as the Navier friction boundary condition,
 \begin{equation} \derp{y}u^\eps_1|_{y=0} = a^\eps u^\eps_1|_{y=0} . \label{NBC} \end{equation}
 This was derived by H. Navier in the XIX\textsuperscript{th} century \cite{Nh} by taking into account the molecular interactions with the boundary.
 To be precise, the Navier condition expresses proportionality between the tangential part of the normal stress tensor and the tangential velocity, thus prescribing how the fluid may slip along the boundary.
 As indicated, the coefficient $a^\eps$ may depend on the viscosity. Typically, we will look at
 \begin{equation} a^\eps = \frac{a}{\eps^\beta}, \label{abeta} \end{equation}
 with $a>0$ and $\beta\geq 0$. A previous paper by the second author \cite{Pm11} showed that nonlinear instability remains present for this type of boundary condition, in particular for the case of boundary-layer-scale data,
 $\beta=1/2$, where there is strong nonlinear instability in $L^\infty$ in the inviscid limit. However, the same article also showed general convergence in $L^2$ when $\beta<1$.
 
\begin{theo}[Theorem 1.2 in \cite{Pm11}]
 \label{theostab} Let $u^\eps_0\in L^2(\Omega)$ and $u^\eps$ be the Leray solution of (\ref{NS}) with initial data $u^\eps_0$, satisfying the Navier boundary conditions (\ref{EBC}) and (\ref{NBC}), with $a^\eps$ as in (\ref{abeta}) with $\beta<1$.
  Let $v_0\in H^s(\Omega)$ with $s>2$, so that $v$ is a global strong solution of the Euler equation (\ref{IE})-(\ref{EBC}), and assume that $u^\eps_0$ converges to $v_0$ in $L^2(\Omega)$ as $\eps\rightarrow0$.
  Then, for any $T>0$, we have the following convergence result:
  $$ \sup_{t\in[0,T]} \norme{L^2(\Omega)}{u^\eps(t)-v(t)} = {\cal O}(\eps^{(1-\beta)/2}). $$
\end{theo}

This theorem is proved using elementary energy estimates and Gr\"onwall's lemma, and it extended results by D. Iftimie and G. Planas \cite{IP}, and X-P. Wang, Y-G. Wang and Z. Xin \cite{WWX}.
 It is worth noting, on one hand, that convergence breaks down for $\beta=1$, and on the other, that a comparable result is impossible to achieve in the no-slip case, since the boundary term $\int_{\partial\Omega} \derp{y}u^\eps_1 v_1~dx$ cannot be dealt with.
\vspace{12pt}

The first remark is important since $\beta=1$ is what we call the ``physical'' case, because this was the dependence on the viscosity predicted by Navier in \cite{Nh}, and because it is indeed
 the Navier condition that one obtains when deriving from kinetic models with a certain scaling (see \cite{MSR} for the Stokes-Fourier system, and recently \cite{JM16} extended the result to Navier-Stokes-Fourier).
 One purpose of this work is therefore to further explore whether or not $\beta=1$ is effectively critical for convergence.
 By using the $L^2$ convergence rate and interpolation, we can obtain a range of numbers $p$ for which convergence in $L^p(\Omega)$ occurs depending on $\beta$, which also breaks down when $\beta=1$. The following extends Theorem \ref{theostab}.

\begin{theo}
 \label{lpconv} Let $u^\eps_0\in L^2(\Omega)$ and $u^\eps$ be the Leray solution of (\ref{NS}) with initial data $u^\eps_0$, satisfying the Navier boundary conditions (\ref{EBC}) and (\ref{NBC}), with $a^\eps$ as in (\ref{abeta}) with $\beta<1$.
  Let $v_0\in H^s(\Omega)$ with $s>2$, so that $v$ is a global strong solution of the Euler equation (\ref{IE})-(\ref{EBC}), and assume that $u^\eps_0$ converges to $v_0$ in $L^2(\Omega)$ as $\eps\rightarrow0$.
  Then, for any $T>0$, we have the following convergence result:
  $$ \lim_{\eps\rightarrow 0} \sup_{t\in[0,T]} \norme{L^p(\Omega)}{u^\eps(t)-v(t)} = 0 ~~~~if~~2\leq p < \frac{2(1+3\beta)}{5\beta-1}. $$
  The convergence rate is $\eps^{(1-\beta)/2-(p-2)(1+3\beta)/4p}$.
\end{theo}

On the second remark, relating to the Dirichlet case, even if no general result like Theorem \ref{theostab} is known, there are necessary and sufficient criteria for $L^2$ convergence. We sum two of these up in the following statement.
\begin{theo}
 \label{katomatsui} Let $u^\eps_0\in L^2(\Omega)$ and $u^\eps$ be the Leray solution of (\ref{NS}) with initial data $u^\eps_0$, satisfying the Dirichlet boundary condition (\ref{DBC}).
  Let $v_0\in H^s(\Omega)$ with $s>2$, so that $v$ is a global strong solution of the Euler equation (\ref{IE})-(\ref{EBC}), and assume that $u^\eps_0$ converges to $v_0$ in $L^2(\Omega)$ as $\eps\rightarrow0$.
  Then, for any $T>0$, the following propositions are equivalent:
  \begin{enumerate}[\bf a.]
   \item $\displaystyle \lim_{\eps\rightarrow 0} \sup_{t\in[0,T]} \norme{L^2(\Omega)}{u^\eps(t)-v(t)} = 0$;
   \item $\displaystyle \lim_{\eps\rightarrow 0} \sqrt{\eps}\int_0^T \norme{L^2(\Gamma_{\kappa\eps})}{\partial_y u_1^\eps(t)}~dt = 0$, where $\Gamma_{\kappa\eps}=\{(x,y)\in\Omega~|~y<\kappa\eps\}$ for $\kappa$ smaller than some $\kappa_0\leq 1$ (a variant of T. Kato \cite{Kt});
   \item $\displaystyle \lim_{\eps\rightarrow 0} \eps \int_0^T \int_{\partial\Omega} (v_1 \derp{y}u^\eps_1)|_{y=0} ~dx~dt = 0$ (S. Matsui \cite{Ms}, Theorem 3).
  \end{enumerate}
\end{theo}
Regarding the key statement {\bf b.} in Theorem \ref{katomatsui}, the original condition found by Kato \cite{Kt} was 
\begin{align}
\lim_{\eps \rightarrow0} \eps \int_0^T \| \nabla u^\eps (t) \|_{L^2 (\Gamma_{\kappa\eps})}^2 d t =0. 
\end{align}
This criterion has been refined by several authors: R. Temam and X. Wang \cite{TW}, X. Wang \cite{W}, J. P. Kelliher \cite{Ke}, and P. Constantin, I. Kukavica, and V. Vicol \cite{CKV}. In fact, the argument of \cite{Kt} provides the inequality 
\begin{align}
\begin{split}
& \limsup_{\eps\rightarrow 0} \sup_{t\in [0,T]} \| u^\eps (t) - v(t) \|_{L^2(\Omega)}^2 \\
& \quad \leq C e^{2 \int_0^T \| \nabla v \|_{L^\infty (\Omega)} d t}  \limsup_{\eps\rightarrow 0}  \eps \left| \int_0^T \langle \partial_y u_1^\eps, {\rm rot}\, \tilde V^{\kappa \eps} \rangle _{L^2(\Omega)}  ~ d t \right|. \label{Kt.ineq}
\end{split}
\end{align} 
Here $C$ is a numerical constant and $\tilde V^{\kappa\eps} (t,x,y) = \tilde V (t,x,\frac{y}{\kappa \eps})$, with a sufficiently small $\kappa\in (0,1]$, is the boundary layer corrector used in \cite{Kt}.
Indeed, Kato's result relied on the construction of a boundary layer at a different scale than in the ansatz presented earlier. It involved an expansion like this,
 $$ u^\eps(t,x,y) = v(t,x,y) + \tilde{V}\left(t,x,\frac{y}{\kappa \eps}\right), $$
 thus convergence in the Dirichlet case is governed by the vorticity's behaviour in a much thinner layer than the physical boundary layer.
 The direction from {\bf b.} to {\bf a.} follows from \eqref{Kt.ineq}.
 Meanwhile, Matsui's result is proved using the energy estimates. 
 
We will show that Theorem \ref{katomatsui} extends `as is' to the Navier boundary condition case.
\begin{theo}
 \label{kmextended} Let $u^\eps_0\in L^2(\Omega)$ and $u^\eps$ be the Leray solution of (\ref{NS}) with initial data $u^\eps_0$, satisfying the Navier boundary conditions (\ref{EBC}) and (\ref{NBC}) with $a^\eps\geq 0$.
  Let $v_0\in H^s(\Omega)$ with $s>2$, so that $v$ is a global strong solution of the Euler equation (\ref{IE})-(\ref{EBC}), and assume that $u^\eps_0$ converges to $v_0$ in $L^2(\Omega)$ as $\eps\rightarrow0$.
  Then, for any $T>0$, convergence in $L^\infty(0,T;L^2(\Omega))$ as in Theorem \ref{theostab} is equivalent to the same Kato and Matsui criteria in the sense as in Theorem \ref{katomatsui}.
\end{theo}

Indeed, we will show that \eqref{Kt.ineq} is valid also for the case of Navier boundary conditions \eqref{EBC} and \eqref{NBC}.
Since the right-hand side of \eqref{Kt.ineq} is bounded from above by 
\begin{align*}
& C e^{2 \int_0^T \| \nabla v \|_{L^\infty (\Omega)} d t}  \limsup_{\eps\rightarrow 0} \kappa^{-\frac12} \eps^\frac12 \int_0^T \| \partial_y u_1^\eps \|_{L^2(\Omega)} d t \\
& \quad \leq C e^{2 \int_0^T \| \nabla v \|_{L^\infty (\Omega)} d t}  \kappa^{-\frac12} \limsup_{\eps\rightarrow 0} \| u^\eps_0 \|_{L^2(\Omega)} T^\frac12. 
\end{align*}
As a direct consequence, we have 
\begin{cor}\label{cor.katomatsui} Under the assumptions of Theorem \ref{katomatsui} or \ref{kmextended}, we have
\begin{align}
\begin{split}
\limsup_{\eps\rightarrow 0} \sup_{t\in [0,T]} \| u^\eps (t) - v(t) \|_{L^2(\Omega)} & \leq C e^{ \int_0^T \| \nabla v \|_{L^\infty (\Omega)} d t}  \| v_0 \|_{L^2(\Omega)}^\frac12 T^\frac14, \label{est.cor}
\end{split}
\end{align}
for some numerical constant $C$.
\end{cor}
Estimate \eqref{est.cor} shows that the permutation of limits 
$$ \lim_{T\rightarrow 0} \lim_{\eps\rightarrow 0} \sup_{t\in [0,T]} \| u^\eps (t) - v(t) \|_{L^2(\Omega)} =\lim_{\eps\rightarrow 0} \lim_{T\rightarrow 0} \sup_{t\in [0,T]} \| u^\eps (t) - v(t) \|_{L^2(\Omega)} $$
is justified, and that this limit is zero, which is nontrivial since $\eps\rightarrow 0$ is a singular limit.
 In particular, at least for a short time period but independent of $\eps$, the large part of the energy of $u^\eps (t)$ is given by the Euler flow $v(t)$.  
\vspace{12pt}

Initially, we hoped to get a result with a correcting layer which could be more tailor-made to fit the boundary condition, but it appears that Kato's Dirichlet corrector yields the strongest statement.
 Whenever we change the $\eps$-scale layer's behaviour at the boundary, we end up having to assume both Kato's criterion and another at the boundary.
 So this result is actually proved identically to Kato's original theorem, and we will explain why in section 3. We will also see that Matsui's criterion extends with no difficulty, but it has more readily available implications.
 
Indeed, the Navier boundary condition gives information on the value of $\derp{y}u^\eps_1$ at the boundary. Assuming that $a^\eps = a\eps^{-\beta}$ as in (\ref{abeta}), we see that
 $$ \eps (v_1 \derp{y}u^\eps_1)|_{y=0} = \eps^{1-\beta} (v_1 u^\eps_1)|_{y=0}. $$
 Simply applying the Cauchy-Schwarz inequality to the integral in the Matsui criterion and using the energy inequality of the Euler equation, we have
 \begin{equation} \eps \int_0^T \int_{\partial\Omega} (v_1 \derp{y}u^\eps_1)|_{y=0}~dx~dt \leq \eps^{1-\beta} \norme{L^2(\Omega)}{v_0} \int_0^T \norme{L^2(\partial\Omega)}{u^\eps_1(t)|_{y=0}}~dt. \label{matnav} \end{equation}
 As the energy inequality for Leray solutions of the Navier-Stokes equation with the Navier boundary condition shows that
 $$ \eps^{1-\beta} \int_0^T \norme{L^2(\partial\Omega)}{u^\eps_1(t)}^2~dt \leq \norme{L^2(\Omega)}{u^\eps(0)}^2, $$
 the right-hand side of (\ref{matnav}) behaves like $C\eps^{(1-\beta)/2}$, and thus converges to zero when $\beta<1$. The Matsui criterion therefore confirms Theorem \ref{theostab}, without being able to extend it to the physical case.
 Once again, the physical slip rate appears to be critical.

\section{Proof of $L^p$ convergence}

To prove Theorem \ref{lpconv}, we rely on \textit{a priori} estimates in $L^\infty$ and interpolation.

First, since the vorticity, $\omega^\eps = \derp{x}u^\eps_2-\derp{y}u^\eps_1$, satisfies a parabolic transport-diffusion equation, the maximum principle shows that
\begin{equation} \norme{L^\infty((0,T)\times \Omega)}{\omega^\eps} \leq \max (\norme{L^\infty(\Omega)}{\omega^\eps|_{t=0}},a\eps^{-\beta}\norme{L^\infty((0,T)\times\partial\Omega)}{u^\eps_1|_{y=0}}) \label{maxp} \end{equation}
by the Navier boundary condition (\ref{NBC})-(\ref{abeta}). To estimate $u^\eps_1$ on the boundary, we use the Biot-Savart law:
$$ u^\eps_1(t,x,0) = \frac{1}{2\pi}\int_{\Omega} \frac{y'}{|x-x'|^2+|y'|^2} \omega^\eps(t,x',y')~dx'dy'. $$
Let us denote $\kappa(x,x',y')$ the kernel in this formula. We split the integral on $y'$ into two parts, $\int_0^K$ and $\int_K^{+\infty}$ with $K$ to be chosen. On one hand, we have
\begin{equation} \left|\int_0^K \int_{\mathbb{R}} \frac{y'}{|x-x'|^2+|y'|^2}\omega^\eps(t,x',y')~dx'dy'\right| \leq C_0 K \norme{L^\infty((0,T)\times\Omega)}{\omega^\eps} \label{omegainf} \end{equation}
by integrating in the variable $x'$ first and recognising the derivative of the arctangent function.

On the other, we integrate by parts, integrating the vorticity $\omega^\eps$, so
\begin{align*}
\begin{split}
& \int_K^{+\infty} \int_{\mathbb{R}} \kappa(x,x',y') \omega^\eps(t,x',y')~dx'dy' \\
& \quad = -\int_K^{+\infty} \int_{\mathbb{R}} u^\eps \cdot \nabla_{x',y'}^\bot \kappa~dx'dy' + \int_{\mathbb{R}} (\kappa u^\eps_1)|_{y'=K}~dx'. 
\end{split}
\end{align*}
The first two terms are easily controlled using the Cauchy-Schwarz inequality: $\norme{L^2}{u^\eps(t)}$ is uniformly bounded by the energy estimate for weak solutions of Navier-Stokes,
 while quick explicit computations show that $\norme{L^2}{\nabla_{x',y'} \kappa}\leq C/K$. Likewise, in the boundary term, the kernel is also ${\cal O}(1/K)$ in $L^2(\mathbb{R})$,
 but we must now control the $L^2$ norm of the trace of $u^\eps_1$ on the set $\{y'=K\}$: by the trace theorem and interpolation, we have
$$ \norme{L^2(\{y'=K\})}{u^\eps_1} \leq \sqrt{\norme{L^2(\Omega)}{u^\eps_1}\norme{L^2(\Omega)}{\omega^\eps}}, $$
and both of these are uniformly bounded. Hence, in total,
$$ \norme{L^\infty((0,T)\times\Omega)}{\omega^\eps} \leq \norme{L^\infty(\Omega)}{\omega^\eps(0)} + a \eps^{-\beta} C_0 K \norme{L^\infty((0,T)\times\Omega)}{\omega^\eps} + a\eps^{-\beta} \frac{C}{K}. $$
By choosing $K \sim \eps^{\beta}$ so that $a\eps^{-\beta}C_0 K < \frac{1}{2}$, we can move the second term on the right-hand side to the left, and we conclude that, essentially,
$$ \norme{L^\infty((0,T)\times\Omega)}{\omega^\eps} \leq C\eps^{-2\beta}. $$

Using the Gagliardo-Nirenberg interpolation inequality from \cite{Nl}, we can now write that, for $p\geq 2$,
$$ \norme{L^p(\Omega)}{u^\eps(t)-v(t)} \leq C\norme{L^2}{u^\eps(t)-v(t)}^{1-q} \norme{L^\infty}{\rot (u^\eps-v)(t)}^q, $$
where $q=\frac{p-2}{2p}$. By Theorem \ref{theostab}, the first term of this product converges to zero with a rate $\eps^{(1-q)(1-\beta)/2}$ when $\beta<1$, while we have just shown that the second behaves like $\eps^{-2q\beta}$, so the bound is
$$ \norme{L^p(\Omega)}{u^\eps(t)-v(t)} \leq C\eps^{(1-\beta)/2 - q(1+3\beta)/2} . $$

It remains to translate this into a range of numbers $p$ such that this quantity converges, which happens when $q<\frac{1-\beta}{1+3\beta}$. Recalling the value of $q$, we get that weak solutions of the Navier-Stokes equation
 converge in $L^p$ towards a strong solution of the Euler equation if
 $$ 2 \leq p < \frac{2(1+3\beta)}{5\beta-1}, $$
 and the right-hand bound is equal to 2 when $\beta=1$.

\section{About the Kato and Matsui criteria}

The starting point for both criteria is the weak formulation for solutions of the Navier-Stokes equation.

\begin{nota} If $E$ is a function space on $\Omega$, we denote $E_\sigma$ the set of 2D vector-valued functions in $E$ that are divergence free and tangent to the boundary.
 Recall that, through the rest of the paper, $a^\eps$ is a non-negative function of $\eps>0$ (not necessarily the same form as in (\ref{abeta})). \end{nota}

\begin{defi}
 A vector field $u^\eps: [0,T]\times\Omega\rightarrow \mathbb{R}^{2}$ is a Leray solution of the Navier-Stokes equation (\ref{NS}) with Navier boundary conditions (\ref{EBC})-(\ref{NBC}) if:
\begin{enumerate}
\item $u^\eps \in {\cal C}_{w}([0,T],L^{2}_\sigma) \cap L^{2}([0,T],H^{1}_\sigma)$ for every $T>0$,
\item for every $\varphi\in H^{1}([0,T],H^{1}_{\sigma})$, we have
\begin{align} \begin{split}
 & \langle u^\eps(T),\varphi(T) \rangle_{L^2(\Omega)}-\int_0^T \langle u^\eps , \partial_{t}\varphi \rangle_{L^2(\Omega)} + \eps a^\eps\int_0^T \int_{\partial\Omega} (u^\eps_1 \varphi_1)|_{y=0} \\
 & \quad + \eps\int_0^T \langle \omega^\eps , \rot \varphi \rangle_{L^2(\Omega)} - \int_0^T \langle u^\eps \otimes u^\eps , \nabla \varphi \rangle_{L^2(\Omega)} = \langle u^\eps(0) , \varphi(0) \rangle_{L^2(\Omega)}, \label{WF} \end{split} \end{align}
\item and, for every $t\geq 0$, $u^\eps$ satisfies the following energy equality (in 3D, this is an inequality):
\begin{equation} \frac{1}{2}\norme{L^{2}(\Omega)}{u^\eps(t)}^{2}+\eps a^{\eps}\int_{0}^{t}\int_{\partial\Omega} (|u^\eps_1|^{2})|_{y=0}+\eps\int_{0}^{t} \norme{L^2(\Omega)}{\omega^\eps}^{2} = \frac{1}{2}\norme{L^{2}(\Omega)}{u^\eps(0)}^{2} . \label{EE} \end{equation}
\end{enumerate}
\end{defi}

When formally establishing the weak formulation (\ref{WF}), recall that
$$ -\int_{\Omega} \Delta u^\eps \varphi = \int_{\Omega} (\omega^\eps \cdot \rot \varphi - \nabla \div u^\eps \cdot \varphi) + \int_{\partial\Omega} (\omega^\eps \varphi\cdot n^\bot)|_{y=0}, $$
where $n^\bot = (n_2, -n_1)$ is orthogonal to the normal vector $n$. In the flat boundary case with condition (\ref{NBC}) on the boundary, we get the third term of (\ref{WF}).
 The differences with the Dirichlet case are two-fold: first, the class of test functions is wider (in the Dirichlet case, the test functions must vanish on the boundary), and second, there is a boundary
 integral in (\ref{WF}) and (\ref{EE}) due to $u^\eps_1$ not vanishing there.
\vspace{12pt}
 
We will not go into great detail for the proof of Theorem \ref{kmextended}, since it is virtually identical to Theorem \ref{katomatsui}. In particular, Matsui's criterion is shown with no difficulty,
 as only the boundary term in (\ref{WF}), with $\varphi=u^\eps-v$, is added in the estimates, and this is controlled as a part of the integral $I_3$ in equality (4.2) in \cite{Ms}, page 167.
 This proves the equivalence \textbf{a.}$\Leftrightarrow$\textbf{c.}
\vspace{12pt}
 
We take more time to show the equivalence \textbf{a.}$\Leftrightarrow$\textbf{b.}, Kato's criterion. In \cite{Kt}, Kato constructed a divergence-free corrector $\tilde V^{\kappa\eps}$, acting at a range ${\cal O}(\eps)$
 of the boundary and such that $v|_{y=0} = \tilde V^{\kappa\eps}|_{y=0}$, and used $\varphi=v-\tilde V^{\kappa\eps}$ as a test function in (\ref{WF}) to get the desired result. We re-run this procedure, which finally leads to the identity 
\begin{align}\label{Key.eq.1}
\begin{split}
& \langle u^\eps (t), v(t) - \tilde V^{\kappa \eps} (t) \rangle _{L^2(\Omega)}\\
& \quad = \langle u^\eps_0, v_0 \rangle_{L^2(\Omega)} - \langle u^\eps_0, \tilde V^{\kappa \eps} (0)\rangle _{L^2(\Omega)} - \int_0^t \langle u^\eps, \partial_t \tilde V^{\kappa\eps} \rangle _{L^2(\Omega)} \\
& \quad \quad + \int_0^t \langle u^\eps - v, (u^\eps - v) \cdot \nabla v \rangle _{L^2(\Omega)} -\eps \int_0^t \langle \omega^\eps, {\rm rot}\, v\rangle _{L^2 (\Omega)} \\
& \quad \quad  + \eps \int_0^t \langle \omega^\eps, {\rm rot}\, \tilde V^{\kappa \eps} \rangle _{L^2(\Omega)} - \int_0^t \langle u^\eps \otimes u^\eps, \nabla \tilde V^{\kappa \eps} \rangle _{L^2(\Omega)}.
\end{split}
\end{align}
In deriving this identity, one has to use the Euler equations which $v$ satisfies and also $\langle v, (u^\eps - v) \cdot \nabla v\rangle _{L^2(\Omega)} =0$. On the other hand, we have from \eqref{EE},
\begin{align}\label{Key.eq.2}
\| u^\eps (t) - v(t) \|_{L^2(\Omega)}^2 & = \| u^\eps (t) \|_{L^2(\Omega)}^2 + \| v(t) \|_{L^2(\Omega)}^2 - 2\langle u^\eps (t), v(t) - \tilde V^{\kappa \eps} \rangle _{L^2 (\Omega)}  \nonumber \\
& \quad \quad - 2\langle u^\eps (t) , \tilde V^{\kappa \eps} (t) \rangle _{L^2 (\Omega)} \nonumber\\
\begin{split}
& = -2 \eps a^{\eps} \int_0^t \| u^\eps_1 \|_{L^2 (\partial\Omega)}^2 -2 \eps \int_0^t \| \omega^\eps\|_{L^2(\Omega)}^2  \\
& \quad \quad + \| u^\eps_0 \|_{L^2(\Omega)}^2 + \| v_0 \|_{L^2 (\Omega}^2 - 2\langle u^\eps (t) , \tilde V^{\kappa \eps} (t) \rangle _{L^2 (\Omega)}  \\
& \quad \quad - 2\langle u^\eps (t), v(t) - \tilde V^{\kappa \eps} (t) \rangle _{L^2 (\Omega)}
\end{split}
\end{align}
Combining \eqref{Key.eq.1} with \eqref{Key.eq.2}, we arrive at the identity which was essntially used reached by Kato in \cite{Kt} for the no-slip case:
\begin{align}\label{Key.eq.3}
\begin{split}
& \| u^\eps (t) - v(t) \|_{L^2(\Omega)}^2 \\
& \quad = -2 \eps a^\eps \int_0^t \| u^\eps_1 \|_{L^2 (\partial\Omega)}^2 -2 \eps \int_0^t \| \omega^\eps\|_{L^2(\Omega)}^2 + \| u^\eps_0 -v_0 \|_{L^2(\Omega)}^2  \\
& \quad \quad - 2\langle u^\eps (t) , \tilde V^{\kappa \eps} (t) \rangle _{L^2 (\Omega)} + 2\langle u^\eps_0, \tilde V^{\kappa\eps} (0) \rangle _{L^2(\Omega)} \\
& \quad \quad  +2 \int_0^t \langle u^\eps, \partial_t \tilde V^{\kappa\eps} \rangle _{L^2(\Omega)} + 2\eps \int_0^t \langle \omega^\eps, {\rm rot}\, v\rangle _{L^2 (\Omega)}\\
& \quad \quad -2 \int_0^t \langle u^\eps - v, (u^\eps - v) \cdot \nabla v \rangle _{L^2(\Omega)} \\
& \quad \quad  +2 \int_0^t \langle u^\eps \otimes u^\eps, \nabla \tilde V^{\kappa \eps} \rangle _{L^2(\Omega)}-2 \eps \int_0^t \langle \omega^\eps, {\rm rot}\, \tilde V^{\kappa \eps} \rangle _{L^2(\Omega)}.
\end{split}
\end{align}
Let us run down the terms in this equality. The first line is comprised of negative terms and the initial difference, which is assumed to converge to zero.
 The terms on the second and third lines of (\ref{Key.eq.3}) tend to zero as $\eps\rightarrow 0$ with the order $\mathcal{O}((\kappa\eps)^\frac12)$, since the boundary corrector has the thickness $\mathcal{O} (\kappa \eps)$.
 Meanwhile, on the fourth line, we have
\begin{align*}
-2 \int_0^t \langle u^\eps - v, (u^\eps - v) \cdot \nabla v \rangle _{L^2(\Omega)}  \leq 2\int_0^t \|\nabla v\|_{L^\infty} \| u^\eps - v\|_{L^2(\Omega)}^2,
\end{align*}
which will be harmless when we apply the Gr\"onwall inequality later. For the Navier-slip condition case, a little adaptation is necessary to control the fifth line,
$$ {\cal I} := \int_0^t \langle u^\eps \otimes u^\eps, \nabla \tilde V^{\kappa\eps} \rangle_{L^2(\Omega)} . $$
In the Dirichlet case, the nonlinear integral ${\cal I}$ is bounded by using the Hardy inequality, since $u^\eps$ vanishes on the boundary. In our case with the Navier condition, however, $u^\eps_1$ does not vanish, so we need to explain this part.
\vspace{12pt}

Let us first manage the terms in ${\cal I}$ which involve $u^\eps_2$, which does vanish on the boundary. Recall that $\tilde V^{\kappa\eps}$ has the form $\tilde V (t,x,\frac{y}{\kappa\eps})$ and is supported in $\Gamma_{\kappa\eps}=\{(x,y)\in\Omega~|~0<y<\kappa \eps\}$, so we write
$$ \left|\int_\Omega (u^\eps_2)^2 \derp{y} \tilde V_2^{\kappa \eps} \right| = \left|\int_{\Gamma_{\kappa\eps}} \left(\frac{u^\eps_2}{y}\right)^2 y^2 \derp{y}\tilde  V_2^{\kappa \eps} \right| \leq C \norme{L^\infty}{y^2 \derp{y} \tilde V_2^{\kappa\eps} } \norme{L^2(\Omega)}{\nabla u^\eps_2}^2, $$
in which we have used the Hardy inequality. Note that $\derp{y} \tilde V^{\kappa\eps}$ is of order $(\kappa \eps)^{-1}$, so $y^2 \derp{y}\tilde V_2^{\kappa\eps}$ is bounded by $C\kappa\eps$ in $L^\infty(\Gamma_{\kappa\eps})$, and we conclude that
\begin{equation} \left|\int_\Omega (u^\eps_2)^2 \derp{y} \tilde V_2^{\kappa\eps}\right| \leq C\kappa \eps \norme{L^2(\Omega)}{\nabla u^\eps}^2. \label{hardy} \end{equation}
Here $C$ is a numerical constant.
 This is what happens on all terms in \cite{Kt}, and the same trick works for $\int_\Omega u^\eps_1 u^\eps_2 \derp{x} \tilde V_2^{\kappa\eps}$;
 this term is in fact better, since the $x$-derivatives do not make us lose uniformity in $\eps$. Using the fact that $\norme{L^2}{u^\eps}$ is bounded courtesy of the energy estimate (\ref{EE}), we have
$$ \left|\int_\Omega u^\eps_1 u^\eps_2 \derp{x} \tilde V_2^{\kappa\eps}\right| \leq C\kappa \eps  \| u^\eps \|_{L^2(\Omega)} \norme{L^2(\Omega)}{\nabla u^\eps}. $$

The term $\int_\Omega u^\eps_1 u^\eps_2 \derp{y} \tilde V_1^{\kappa\eps}$ is trickier, since the $y$-derivative is bad for uniformity in $\eps$, and we only have one occurrence of $u^\eps_2$ to compensate for it.
 Let us integrate this by parts: using the divergence-free nature of $u^\eps$, we quickly get
\begin{eqnarray*} \int_\Omega u^\eps_1 u^\eps_2 \derp{y} \tilde V_1^{\kappa\eps} & = & \int_\Omega u^\eps_1 \derp{x}u^\eps_1 \tilde V_1^{\kappa\eps} - \int_\Omega \derp{y}u^\eps_1 u^\eps_2 \tilde V_1^{\kappa\eps} \\
 & = & -\frac{1}{2} \int_\Omega (u^\eps_1)^2 \derp{x} \tilde V_1^{\kappa\eps} - \int_\Omega \derp{y}u^\eps_1 u^\eps_2 \tilde V_1^{\kappa\eps}. \end{eqnarray*}
The second term can be dealt with using the Hardy inequality as above, and its estimate is identical to (\ref{hardy}). The first term, meanwhile, is the same as the remaining one in ${\cal I}$.

To handle $\int_\Omega (u^\eps_1)^2 \derp{x} \tilde V_1^{\kappa\eps}$, in which no term vanishes on the boundary, we proceed using the Sobolev embedding and interpolation. Indeed, we have
\begin{align}
\left|\int_\Omega (u^\eps_1)^2 \derp{x} \tilde V_1^{\kappa\eps} \right| & \leq 2 \norme{L^2(\Omega)}{(u^\eps_1 - v_1)^2} \norme{L^2(\Omega)}{ \partial_x \tilde V_1^{\kappa\eps}} + 2 \| v_1^2 \|_{L^2(\Omega)} \| \partial_x \tilde V_1^{\kappa\eps}\|_{L^2(\Omega)} \nonumber \\
& \leq C (\kappa \eps)^\frac12 \| u_1^\eps - v_1 \|_{L^4(\Omega)}^2 + C \kappa\eps \| v\|_{L^\infty(\Omega)}^2.
\end{align} 
Here we have used that $\norme{L^2(\Omega)}{ \tilde V_1^{\kappa\eps}} \leq C(\kappa\eps)^\frac12$, while
$$ \norme{L^4(\Omega)}{u^\eps_1-v_1}^2  \leq C\norme{L^2(\Omega)}{u^\eps_1-v_1}\norme{H^1(\Omega)}{u^\eps_1-v_1}, $$
and so, in total, we conclude that
\begin{align}\label{est.I}
\begin{split}
|{\cal I}| & \leq C\bigg ( \kappa \eps \norme{L^2(\Omega)}{\nabla u^\eps}^2 +\kappa \eps  \| u^\eps \|_{L^2(\Omega)} \norme{L^2(\Omega)}{\nabla u^\eps}\\
&  \quad + \| u^\eps - v\|_{L^2 (\Omega)}^2 + (\kappa\eps)^\frac12 \| \nabla v \|_{L^2(\Omega)} \| u^\eps-v\|_{L^2(\Omega)}  + \kappa\eps \| v\|_{L^\infty (\Omega)}^2 \bigg). 
\end{split}
\end{align}
Here $C$ is a numerical constant. 
 Then, by virtue of the identity $\|\omega^\eps \|_{L^2(\Omega)} = \| \nabla u^\eps\|_{L^2(\Omega)}$,
 the term $C\kappa \eps \| \nabla u^\eps\|_{L^2 (\Omega)}^2$ in the right-hand side of \eqref{est.I} can be absorbed by the dissipation in the first line of \eqref{Key.eq.3} if $\kappa>0$ is suifficiently small.
\vspace{12pt}

We come to the final linear term $-2\eps\int_0^t \langle \omega^\eps, {\rm rot}\, \tilde V^{\kappa\eps}\rangle _{L^2 (\Omega)}$ in the fifth line of (\ref{Key.eq.3}). Using $\omega^\eps = \partial_x u^\eps_2 - \partial_y u^\eps_1$, we have from the integration by parts,
\begin{align*}
& -2\eps\int_0^t \langle \omega^\eps, {\rm rot}\, \tilde V^{\kappa\eps}\rangle _{L^2 (\Omega)} \\
& \quad = 2\eps \int_0^t \langle \partial_y u_1^\eps, {\rm rot}\, \tilde V^{\kappa\eps}\rangle _{L^2 (\Omega)}  + 2\eps \int_0^t \langle \frac{u_2^\eps}{y}, y \partial_x {\rm rot}\, \tilde V^{\kappa\eps} \rangle _{L^2(\Omega)} \\
& \quad \leq 2\eps \int_0^t \langle \partial_y u_1^\eps, {\rm rot}\, \tilde V^{\kappa\eps}\rangle _{L^2 (\Omega)}  + C  \kappa^\frac12 \eps^\frac32 \int_0^t \| \nabla u^\eps_2 \|_{L^2(\Omega)}.
\end{align*}
Collecting all these estimates, we get from \eqref{Key.eq.3} that for $0<t\leq T$,
\begin{align}\label{Key.ineq}
\begin{split}
& \| u^\eps (t) - v(t) \|_{L^2(\Omega)}^2 \\
& \quad \leq  - 2 \eps a^\eps \int_0^t \| u^\eps_1 \|_{L^2 (\partial\Omega)}^2 - \eps \int_0^t \| \omega^\eps\|_{L^2(\Omega)}^2  + \| u^\eps_0 -v_0 \|_{L^2(\Omega)}^2 \\
& \quad \quad  + C (\kappa\eps)^\frac12  +  \int_0^t \big (C_0 + 2 \| \nabla v \|_{L^\infty (\Omega)} \big ) \| u^\eps - v \|_{L^2(\Omega)}^2 \\
& \quad \quad + 2 \eps \int_0^t \langle \partial_y u_1^\eps, {\rm rot}\, \tilde V^{\kappa \eps} \rangle _{L^2(\Omega)}.
\end{split}
\end{align}
Here $C$ depends only on $T$, $\| u^\eps_0 \|_{L^2(\Omega)}$, and $\|v_0 \|_{H^s(\Omega)}$, while $C_0$ is a numerical constant. Inequality \eqref{Key.ineq} is valid also for the no-slip (Dirichlet) case; indeed, we can drop the negative term $- 2 \eps a^\eps \int_0^t \| u^\eps_1 \|_{L^2 (\partial\Omega)}^2$.
 By applying the Gr\"onwall inequality and by taking the limit $\eps\rightarrow 0$, we arrive at \eqref{Kt.ineq}.
 This is enough to extend Kato's criterion to the Navier boundary condition case; the rest is identical to Kato's proof in \cite{Kt}.
 We have achieved this result by re-using the Dirichlet corrector because, since the test function $\varphi=v-\tilde V^{\kappa\eps}$ vanishes at $y=0$, the boundary integral in (\ref{WF}) does not contribute.
 This does not feel quite satisfactory. One would have hoped to get criteria by constructing more appropriate correctors, such as one so that the total satisfies the Navier boundary condition, but, as we have just mentioned,
 a boundary integral appears and it is not clear that we can control it. In fact, this boundary term is similar to the one in the Matsui criterion, which, as we have proved, is equivalent to Kato's.
 We observe that when considering a corrector which does not vanish on the boundary, convergence of Navier-Stokes solutions to Euler solutions happens if and only if both \textbf{b.} and \textbf{c.} are satisfied.
 It appears difficult to get refinements of criteria for $L^2$ convergence in the inviscid limit problem according to the boundary condition.
\vspace{12pt}

\textbf{Acknowledgements.} These results were obtained during MP's visit to Kyoto University, the hospitality of which is warmly acknowledged. The visit was supported by the JSPS Program for Advancing Strategic International Networks to Accelerate the Circulation of Talented Researchers, 
 `Development of Concentrated Mathematical Center Linking to Wisdom of the Next Generation', which is organized by the Mathematical Institute of Tohoku University.
 MP is also supported by the European Research Council (ERC) under the European Union's Horizon 2020 research and innovation program Grant agreement No 63765, project `BLOC', as well as the French Agence Nationale de la Recherche project `Dyficolti' ANR-13-BS01-0003-01.

\begin{small}
\bibliography{notebib.10.5}
\bibliographystyle{abbrv}
\end{small}

\end{document}